\documentclass[12pt]{article}

\newcommand\bC{\mathbf{C}}
\newcommand\bP{\mathbf{P}}
\newcommand\bQ{\mathbf{Q}}
\newcommand\bZ{\mathbf{Z}}

\newcommand\cC{\mathcal{C}}
\newcommand\cI{\mathcal{I}}
\newcommand\cJ{\mathcal{J}}
\newcommand\cM{\mathcal{M}}
\newcommand\cO{\mathcal{O}}

\newtheorem{conj}{Conjecture}
\newtheorem{lem}{Lemma}
\newtheorem{prop}{Proposition}

\renewcommand\b{\beta}
\renewcommand\l{\lambda}

\newcommand\ch{\mathrm{ch}}
\newcommand\Ext{\mathrm{Ext}}
\newcommand\mod{\mathrm{mod}}
\newcommand\DT{\mathrm{DT}}

\begin{document}
  \begin{center}
{\LARGE    Gromov-Witten, Gopakumar-Vafa, and Donaldson-Thomas invariants of
Calabi-Yau threefolds}\\
\ \\
Sheldon Katz\\
Departments of Mathematics and Physics\\
University of Illinois at Urbana-Champaign\\
Urbana, IL 61801
  \end{center}

\bigskip\noindent 
{\em Abstract:\/} Gromov-Witten, Gopakumar-Vafa, and
Donaldson-Thomas invariants of Calabi-Yau threefolds are compared.  In
certain situations, the Donaldson-Thomas invariants are very easy to handle,
sometimes easier than the other invariants.
This point is illustrated in several ways,
especially by revisiting computations of Gopakumar-Vafa invariants by
Katz, Klemm, and Vafa in a rigorous mathematical framework.  This note
is based on my talk at the 2004 Snowbird Conference on String
Geometry.

\bigskip

\section{DT Invariants and GW Invariants}
\subsection{Generalities}
Let $X$ be a nonsingular complex projective threefold, $\b\in
H_2(X,\bZ)$, and let $n\in\bZ$.  We let $I_n(X,\b)$ denote the part of
the Hilbert scheme of $X$ parametrizing subschemes $Z\subset X$ with
\begin{itemize}
  \item $[Z]=\b$
\item $\chi(\cO_Z)=n$.
\end{itemize}
The class $[Z]\in H^2(X,\bZ)$ can be equivalently defined as either
the dimension one component of the support cycle of $Z$, or as
$\ch_2(\cO_Z)$.

We let $\cI_Z$ be the ideal sheaf of $Z$.  In \cite{thomas}, a perfect
obstruction theory is defined on $I_n(X,\b)$ arising naturally from
the deformation theory of the ideal sheaves $\cI_Z$.  The virtual dimension
is given by
\begin{equation}
  \label{vd}
D=\dim\Ext^1_0(\cI_Z,\cI_Z)-\dim\Ext^2_0(\cI_Z,\cI_Z)=c_1(X)\cdot\b.
\end{equation}
In (\ref{vd}), the subscript of 0 denotes traceless Ext, i.e.\ the
kernel of the trace map:
\begin{equation}
  \label{trext}
\Ext^i_0(\cI_Z,\cI_Z)=\ker\left(
\Ext^i(\cI_Z,\cI_Z)\to H^i(X,\cO)
\right).
\end{equation}

Associated to the perfect obstruction theory, there is a virtual fundamental
class
\begin{equation}
  \label{vfc}
[I_n(X,\b)]^\mathrm{vir}\in H_{2D}\left(I_n(X,\b)
\right).
\end{equation}

If $\b=0$, then the virtual dimension $D$ is 0 for all $n$.  
Putting
\begin{equation}
  \label{d0n}
D_0^n=\deg[I_n(X,0)]^\mathrm{vir}\in\bZ
\end{equation}
and introducing a formal variable $q$,
the {\em dimension~0 Donaldson-Thomas partition function\/} is defined as the
formal series
\begin{equation}
  \label{zdt0}
Z_0^{\mathrm{DT}}(X)=\sum_{n=0}^\infty D_0^nq^n.
\end{equation}

\begin{conj}[\cite{mnop}]
\label{conj1}
  \[
Z_0^{\mathrm{DT}}(X)=M(-q)^{\int_Xc_3-c_1c_2}.
\]
\end{conj}
Here $M(q)$ is the McMahon function, the generating function of 
three dimensional partitions, 
\begin{equation}
  \label{mcmahon}
M(q)=\prod_{n=1}^\infty\left(1-q^n\right)^{-n}.
\end{equation}
The $c_i$ are the Chern classes of $X$.  In particular, $\int_Xc_3=e(X)$,
the topological Euler characteristic of $X$.

Conjecture~\ref{conj1} has been proven if $X$ is a toric variety \cite{mnop2}.
A key feature is the use of localization which relates the conjecture to
a combinatorial counting of 3~dimensional partitions.  The relevance 
of 3~dimensional partitions to Gromov-Witten theory appeared in \cite{orv}.
In the present note, we will explain some simple computations in DT
theory which do not rely on localization.

For each $n$, Conjecture~\ref{conj1} makes a prediction for $D_0^n$.
This prediction is trivially true for $n=1$ and has been proven for $n=2$ by
Maulik and Pandharipande \cite{mp}.  A simpler argument can be given in the
Calabi-Yau case for $n\le3$ as will be described below.

\medskip
For $\b\ne0$, the virtual dimension need not be zero.  Nevertheless, invariants
can be obtained by imposing additional conditions \cite{mnop2} or by
working equivariantly \cite{bp}.  

\medskip
Note that another important way to achieve virtual dimension
$D=0$ for all $n$ (and all $\b$ as well) is if $X$
is Calabi-Yau, which we take in this paper to simply mean that $K_X$
is trivial.
This case is the focus of the present paper, and from now on we assume that
$X$ is Calabi-Yau.  Note that
the conclusion on the virtual dimension also holds more
generally if $c_1(X)$ is torsion.  

\subsection{Calabi-Yau case}
\label{cy}
Assuming that $X$ is Calabi-Yau, we now put
\begin{equation}
  \label{dbn}
D_\b^n=\deg[I_n(X,\b)]^\mathrm{vir}\in\bZ
\end{equation}
and define the {\em degree $\b$ Donaldson-Thomas partition function\/} by
\begin{equation}
  \label{zdtb}
Z_\b^{\mathrm{DT}}(X)=\sum_n D_\b^nq^n.
\end{equation}
We now introduce formal symbols $t^\b$ satisfying $t^\b t^{\b'}=t^{\b+\b'}$
and define the full {\em Donaldson-Thomas partition function\/}
\begin{equation}
  \label{zdt}
Z^{\mathrm{DT}}(X)=\sum_\b Z_\b^{\mathrm{DT}}(X)t^\b
\end{equation}
and the {\em reduced Donaldson-Thomas partition function\/} 
$Z^{\mathrm{DT}}(X)'$ and
{\em reduced degree $\b$ Donaldson-Thomas partition functions\/}
$Z_\b^{\mathrm{DT}}(X)'$ by 
\begin{equation}
  \label{zdtn}
Z^{\mathrm{DT}}(X)'=\frac{Z^{\mathrm{DT}}(X)}{Z^{\mathrm{DT}}_0(X)}
=:\sum_\b Z_\b^{\mathrm{DT}}(X)'t^\b.
\end{equation}

It is then natural to define the {\em reduced Donaldson-Thomas invariants\/}
${D^n_\b}'$ by
\begin{equation}
  \label{normdt}
Z^{\mathrm{DT}}_\b(X)'=\sum_n {D^n_\b}'q^n.
\end{equation}

\medskip
Now recall the Gromov-Witten invariants $N^g_\b=\deg[\overline{M}_{g,0}
(X,\b)]^\mathrm{vir}\in\bQ$.  Form the generating functions
\[
F'_g(X)=\sum_{\b\ne0}N^g_\b t^\b,\qquad F'(X)=\sum_g\l^{2g-2}F'_g,\qquad
Z^{\mathrm{GW}}(X)'=\mathrm{exp}(F'(X)).
\]
The prime is used to emphasize that the constant maps $\b=0$ are not included.

\begin{conj} (\cite{mnop})
\label{conj2}
$Z^{\mathrm{DT}}(X)'=  Z^{\mathrm{GW}}(X)'$ after the change of variables
$q=-e^{i\l}$.
\end{conj}

In \cite{mnop}, localization techniques were used to prove
Conjecture~\ref{conj2} if $X$ is a noncompact toric Calabi-Yau
threefold, under the additional assumption that the topological vertex
prediction \cite{akmv,orv} correctly computes the Gromov-Witten
partition function.

Note that in \cite{mnop}, a compactification of $X$ is used to
make sense of the Donaldson-Thomas invariants of $X$, and then
$Z^{\mathrm{DT}}(X)'$, appropriately defined, is independent of the
choice of compactification.  The invariants can also be defined in these
non-compact cases by working equivariantly as in \cite{bp}.

\medskip
Donaldson-Thomas invariants are often easy to compute in the Calabi-Yau case.
Suppose $I_n(X,\b)$ is smooth.  Then for $\cI\in I_n(X,\b)$, we have
$\Ext^1_0(\cI,\cI)=T_{I_n(X,\b),\cI}$.  Since $\Ext^2_0(\cI,\cI)\simeq
(\Ext^1_0(\cI,\cI))^*$ canonically (up to a global constant), it follows
that the obstruction bundle is $T^*I_n(X,\b)$.  If furthermore we have
a decomposition
\[
I_n(X,\b)=\cup M_i
\]
into connected components, then we have the simple formula
\begin{equation}
  \label{simple}
D_{\b}^n=\sum_i e(T^*M_i)=\sum_i(-1)^{\dim(M_i)}e(M_i).
\end{equation}
where $e(M_i)$ denotes the topological Euler characteristic of $M_i$.

Note that in the Calabi-Yau case, Conjecture~\ref{conj1} takes the simpler
form
\begin{equation}
\label{conj1cy}
  Z_0^{\mathrm{DT}}(X)=M(-q)^{e(X)}.
\end{equation}

Compare (\ref{conj1}) with the following.

\begin{prop}[\cite{cheah}]
\label{cheah}
  $\sum_{n=0}^\infty e(I_n(X,0))q^n=M(q)^{e(X)}$.
\end{prop}

Note that $I_n(X,0)$ is smooth for $n\le3$ and singular for 
$n\ge4$.\footnote{The
singular locus of $I_4(X,0)$ is described in e.g.\ \cite{hilb4}.}
It follows that for $n\le3$
we have $D^n_0=(-1)^{3n}e(I_n(X,0))=(-1)^ne(I_n(X,0))$.  
It follows immediately from
(\ref{simple}) and (\ref{cheah}) that
$Z_0^{\mathrm{DT}}(X)\equiv M(-q)^{e(X)}\ (\mod\ q^4)$.
The same argument shows that (\ref{conj1cy}) is equivalent to
\[
D_0^n=(-1)^ne(I_n(X,0))
\]
for all $n$.

\bigskip
As a simple illustration of the underlying ideas, we explain why
$D_0^3$ is the coefficient of $q^3$ in $M(-q)^{e(X)}$.

\medskip
First note that $I_3(X,0)$ is the Hilbert scheme of degree three zero
dimensional subschemes in $X$, which is well known to be smooth.  Next,
$I_3(X,0)$ can be obtained from $\mathrm{Sym}^3X$ by two blowups.
First blow up along the small diagonal $\Delta_s$, which is isomorphic
to $X$.  Then blow up the proper transform of the big diagonal $\Delta$
to get $I_3(X,0)$.  We have a disjoint union
\[
I_3(X,0)=\left(\mathrm{Sym}^3X-\Delta\right)\cup
\left(\bP^2{\rm \ bundle\ over\ }\Delta-\Delta_s\right)\cup
\left(\bP^5{\rm \ bundle\ over\ }\Delta_s\right)
\]
This leads to
\[
\begin{array}{ccl}
e\left(I_3(X,0)\right)&=&
\frac16e(X^3-\Delta)+3e(X)(\Delta-\Delta_s)+6e(X)(\Delta_s)\\
&=&\frac16(e(X)^3-3e(X)^2+2e(X))+3(e(X)^2-e(X))+6e(X),
\end{array}
\]
which is immediately checked to be the coefficient of $q^3$ in $M(-q)^{e(X)}$.

A similar and easier calculation, blowing up $\mathrm{Sym}^2X$ along the
diagonal shows that
$D_0^2$ is the coefficient of $q^2$ in $M(-q)^{e(X)}$.

These techniques use the motivic property of the Euler characteristic.
If the dimension zero Donaldson-Thomas invariants can be shown to be
appropriately motivic, then Conjecture~\ref{conj1} would follow from
the computation of the dimension zero Donaldson-Thomas invariants of
$\bC^3$ as computed by localization.  

\smallskip
A more general recipe for obtaining the component of $I_n(X,0)$ containing
distinct points by blowing up a symmetric product of $X$ is given in
\cite{es}.

\section{GV invariants}
From M-theory compactified on $X$ we get conjectured Gopakumar-Vafa
invariants $n^g_\b\in\bZ$ \cite{gv} which satisfy
\begin{equation}
\label{gv}
F'(X)=\sum_{m,g,\b\ne0}
n^g_{\b}\frac1m\left(2\sin\frac{m\lambda}2\right)^{2g-2}t^{m\b}.
\end{equation}
At present, there is no intrinsic mathematical definition of the $n^g_{\b}$.
One can take (\ref{gv}) as a recursive definition of the Gopakumar-Vafa
invariants in terms of the Gromov-Witten invariants, but then it is only
clear that $n^g_{\b}\in\bQ$.  The integrality of the $n^g_\b$ is called the
{\em integrality conjecture\/}.  A discussion from the mathematical
perspective appears in \cite{integrality}.  One can also attempt to 
define the GV invariants recursively using the DT invariants, as will be
illustrated below.  The equivalence of these recursive definitions would 
follow from the validity of Conjecture~\ref{conj2}.

As is well known, Gromov-Witten invariants can be defined on
noncompact Calabi-Yau threefolds in a number of interesting contexts,
e.g.\ neighborhoods of super-rigid curves \cite{superrigid} or
homology classes contained in Fano surfaces.  

\smallskip\noindent {\bf Examples.} Let $X$ be the total space of the
bundle $\cO(-1)\oplus \cO(-1)$ over $\bP^1$ (``local $\bP^1$'').
Identifying $\bP^1$ with the zero section we have $n^0_{[\bP^1]}=1$
and all other $n^g_\b=0$. For a general local elliptic
curve we have $n^1_{k[E]}=1$ for all $k\ge1$, with all other
$n^g_\b=0$.  See \cite{superrigid} for more detail.

\bigskip
Conjecture~\ref{conj2} together with (\ref{gv}) implies that the coefficient
\begin{equation}
\label{coeff}
\sum_m\frac1m\left(2\sin\frac{m\lambda}2
\right)^{2g-2}t^{m\b}
\end{equation} 
of a
genus $g>0$ GV invariant in $F'$ contributes to the 
DT partition function as
\begin{equation}
\begin{array}{ccl}
{Z^{\mathrm{DT}}}'&=&\mathrm{exp}\left(\sum_m\frac1m\left(2\sin\frac{m\lambda}2
\right)^{2g-2}t^{m\b}\right)\\
&=&\mathrm{exp}\left(\sum_m\frac1m\left(-1\right)^{g-1}\left(e^{im\lambda/2}
-e^{-im\lambda/2}
\right)^{2g-2}t^{m\b}\right)\\
&=&\mathrm{exp}\left(\sum_{m,k}\frac1m(-1)^{k+g-1}{2g-2 \choose k}
e^{i(g-1-k)m\lambda}t^{m\b}\right)\\
&=&\mathrm{exp}\left(\sum_k(-1)^{k+g}{2g-2 \choose k}
\log\left(1-e^{i(g-1-k)\lambda}t^\b\right)\right)\\
&=&\prod_{k=0}^{2g-2}\left(1-e^{i\left(g-1-k\right)\lambda}t^\b
\right)^{\left(-1\right)^{k+g}{2g-2 \choose k}}\\
&=&\prod_{k=0}^{2g-2}\left(1+\left(-1\right)^{g-k}q^{g-1-k}t^\b
\right)^{\left(-1\right)^{k+g}{2g-2 \choose k}.}
\end{array}
\end{equation}

If $g=0$, we get 
\begin{equation}
\begin{array}{ccl}
{Z^{\mathrm{DT}}}'&=&\mathrm{exp}\left(\sum\frac1m\left(2\sin\frac{m\lambda}2
\right)^{-2}t^{m\b}\right)\\
&=&\mathrm{exp}\left(\sum_m
-e^{im\lambda}\frac1m \left(1-e^{im\lambda}\right)^{-2}
t^{m\b}\right)\\
&=&\mathrm{exp}\left(\sum_{m,k}-\frac1mke^{imk\lambda}t^{m\b}\right)\\
&=&\prod_{k=1}^\infty\left(1-e^{ik\lambda}t^\b\right)^k\\
&=&\prod_{k=1}^\infty\left(1+(-1)^{k+1}q^kt^\b\right)^k.
\end{array}
\end{equation}

Putting these two cases together, we get for $Z^{\DT}(X)'$ the product 
expression
\begin{equation}
\label{gvdt}
\prod_\b\left(
\prod_{j=1}^\infty\left(1+(-1)^{j+1}q^jt^\b\right)^{jn^0_\b}
\prod_{g=1}^\infty\prod_{k=0}^{2g-2}\left(1+\left(-1\right)^{g-k}q^{g-1-k}t^\b
\right)^{\left(-1\right)^{k+g}n^g_\b{2g-2 \choose k}}\right)
\end{equation}
in terms of the GV invariants.
We get $Z^{\DT}(X)$ by multiplying (\ref{gvdt}) by $M(-q)^{e(X)}$.

\medskip
The factor corresponding to $\b$ in (\ref{gvdt}) can be expanded to
\begin{equation}
  \label{gvdtexp}
1+t^\b\left(n^g_\b q^{1-g}+\left(n^{g-1}_\b+(2g-2)n^g_\b
\right)
+O(q^{3-g})\right)+O(t^{2\b}),
\end{equation}
where $g$ is now interpreted as 
the maximum $g$ such that $n^g_\b$ is nonzero.

\bigskip\noindent
{\bf Example.} The adaptation of 
\cite[Theorem~3]{mnop} to local $\bP^1$ combined with the computation
of the Gromov-Witten invariants of local $\bP^1$ in \cite{fp}
gives 
\[
Z_{\DT}'=\prod_{k=1}^\infty\left(1+(-1)^{k+1}q^kt^\b\right)^k
\]
as predicted.  The Donaldson-Thomas invariants of local $\bP^1$ can also be
computed by the methods of \cite{orv}.

\medskip\noindent {\bf Example.} If $g=1$, then Conjecture~\ref{conj2}
implies that
$Z'_{DT}=\prod_m(1-t^{m\b})^{-1}= \sum p(k) t^{k\b}$, where $p$
is the classical partition function.  We partially verify this prediction
by computing ${D^n_{k[E]}}'$ for $n=0$ and all $k$,
as well as for $n=k=1$, where $E\subset X$ is the elliptic curve.

It is straightforward to reduce to the case where
$X$ is locally $L\oplus L^{-1}$, with $L$ a nontorsion degree zero line
bundle on the elliptic curve $E$.  Since $e(X)=0$, it follows
that the DT invariants and reduced DT invariants coincide in the local case.
Note that $L$ and $L^{-1}$
are the only quotient bundles of $\cI_Z/\cI_Z^2$ of degree 0.  Let
$x$ and $y$ denote the fiber coordinates on $L$ and $L^{-1}$ respectively,
so that in particular the zero section $E$ has equations $x=y=0$.  We can
associate a monomial ideal in $(x,y)$ to any partition of $k$ in the usual
way.  In the current context, this monomial ideal defines a subscheme
$Z\subset X$ supported on $E$, with multiplicity $k$ and $\chi(\cO_Z)=0$.
It can be checked that these are the only such subschemes.  Each of
these ideals 
defines an isolated reduced point of $I_0(X,k[E])$.  Since there are $p(k)$
such points, this contributes $p(k)t^{k\b}$ to the DT partition function.

\medskip
Now, returning to a more general $X$, suppose for simplicity that $E$
is the only curve in $X$ of class $\b$.  If $p\in X-E$, then
$Z=E\cup\{p\}$ has ideal sheaf $\cI_Z\in I_1(X,\b)$.  There are also
ideal sheaves $\cI\in I_1(X,\b)$ with support $E$ and a single embedded
point.  It is not hard to see from this that $I_1(X,\b)=
\mathrm{Bl}_{E}X$ (the proof is much easier than the proof of
Lemma~\ref{blowup} below).  Clearly $e(\mathrm{Bl}_{E}X)= e(X)$, so
$D_\b^1=-e(X)$.  This implies
\[
Z_\b^{DT}(X)\equiv 1-e(X)q\ \mod\ q^2
\]
which implies
\[
Z_\b^{DT}(X)'\equiv 1\ \mod\ q^2,
\]
consistent with the prediction  $Z_\b^{\mathrm{DT}}(X)'=1$.

\bigskip
In the context of DT invariants, we now revisit the
methods for computing GV invariants which were developed in
the physics literature \cite{kkv}.

\smallskip
Fix $g\ge0$, and suppose $I_{1-g}(X,\b)$ parametrizes ideals of local
complete intersection curves of arithmetic genus $g$.\footnote{In
\cite{kkv}, the curves were contained in smooth Fano surfaces
embedded in Calabi-Yau threefolds, e.g.\ ``local $\bP^2$'', but
Proposition~\ref{kkvprop} below applies in greater generality.}  To
compare with \cite{kkv}, we put $\cM=I_{1-g}(X,\b)$.  Identifying
$\cM$ with a component of the Hilbert scheme, let $\cC\subset
\cM\times X$ be the universal subscheme, and let $\cC^{[n]}$ be the
relative Hilbert scheme of $n$ points in the family $\cC/\cM$.  In
particular $\cC^{[1]}= \cC$.  Put $\cC^{[0]}=\cM$.

Fix $\delta\le g$, and assume that $\cC^{[n]}$ is smooth for $n\le\delta$.
Then under certain additional hypotheses

\begin{conj}  (\cite{kkv})
\label{kkv}
\[
\left(-1\right)^{\dim\cM+\delta}n^{g-\delta}_\b=e\left(\cC^{[\delta]}\right)+
\left(2g-2\delta\right)e\left(\cC^{[\delta-1]}\right)+\]
\[\sum_{i=2}^\delta\frac1{i!}\left(2g-2\delta+2i-2\right)
\left(2g-2\delta+i-3\right)\left(2g-2\delta+i-4\right)\cdots
\left(2g-2\delta-1\right)e\left(\cC^{[\delta-i]}\right).
\]
\end{conj}

Using (\ref{gvdtexp}), Conjectures~\ref{conj2} and \ref{kkv} imply that 
\begin{equation}
\label{d1gb}
{D^{1-g}_\b}'=n^g_\b=(-1)^{\dim\cM}e(\cM)
\end{equation}
\begin{equation}
  \label{d2gb}
{D^{2-g}_\b}'=n^{g-1}_\b+(2g-2)e(\cM)=(-1)^{\dim\cC}e(\cC).
\end{equation}

\begin{prop}
\label{kkvprop}
In addition to the above hypotheses, suppose that all ideals $\cI\in
I_{2-g}(X,\b)$ have $\cO_X/\cI$ supported on a curve parametrized by
$\cM$.  Then, for $\delta\le1$ the contribution of $\cM$ to
$Z^{\mathrm{DT}}_\b(X)'$ is
\[
\sum_{n=0}^\delta\left(-1\right)^{\dim\cC^{[n]}}e(\cC^{[n]})q^{n+1-g}+
O(q^{\delta+2-g}).
\]
\end{prop}

Proposition~\ref{kkvprop} has been written in this suggestive way as
omitted computational evidence suggests that Conjecture~\ref{kkv} may
imply the validity of the assertion of Proposition~\ref{kkvprop} with
appropriate additional hypotheses for
larger values of $\delta$.  However, new techniques will be needed to
investigate this situation, as the DT moduli spaces will not be smooth
in general.

The examples in \cite{kkv} satisfy the hypotheses of Proposition~\ref{kkvprop},
which may therefore be viewed as providing a mathematically
rigorous framework for the computational techniques developed in \cite{kkv}.

\begin{lem}
\label{blowup}
  Under the above hypotheses, $I_{2-g}(X,\b)$ is the blowup of
$\cM\times X$ along $\cC$.
In particular, $I_{2-g}(X,\b)$ is smooth.
\end{lem}

\medskip\noindent
{\em Proof:\/} Let $\cI\subset\cO_{\cM\times X}$ denote the universal
ideal sheaf.  Let $\pi_{13}:\cM\times X\times X\to \cM\times X$ be the
projection onto the first and third factors and consider the composition
of mappings
\begin{equation}
\label{compos}
  \pi^*_{13}\cI\hookrightarrow\cO_{\cM\times X\times X}\to \cO_{\cM\times
\Delta_X}
\end{equation}
with $\Delta_X\subset X\times X$ the diagonal and the second map being
the restriction.  Interpreting (\ref{compos}) as a family of maps over
$\cM\times X$, its fiber over $(Z,p)$ is surjective unless
$(Z,p)\in\cC$, so that the ideal of $\cC$ is defined as the
scheme-theoretic locus in $\cM\times X$ over which (\ref{compos})
vanishes.  Let $\pi:\mathrm{Bl}_{\cC}(\cM\times X)\to \cM\times X$ be
the blowup of $\cC$, with $E\subset \mathrm{Bl}_{\cC}(\cM\times X)$ the
exceptional divisor.

Let $\iota:X\times X\to X\times X$ be the map interchanging factors of $X$,
and let $\rho$ denote the composition
\begin{equation}
  \label{rhodef}
\rho:\mathrm{Bl}_{\cC}(\cM\times X)\times X\stackrel{\pi\times
1_X}{\longrightarrow} \cM\times X\times X \stackrel{1_{\cM}\times
\iota}{\longrightarrow} \cM\times X\times X.
\end{equation}
Then (\ref{compos}) induces a mapping
\begin{equation}
  \label{univ}
\phi:\rho^*\cI\to \cO_{\rho^{-1}(\cM\times\Delta_X)}(-E).
\end{equation}
Using local generators of $\cI$ it is straightforward to check that
(\ref{univ}) is surjective.  Locally over $\cM\times X$, let $\cI$ be
generated by $f=f(Z,p,x),\ g=g(Z,p,x)$, where $Z\in\cM$ and $p,x\in X$, so that
$\cC$ is defined by $f(Z,p,p)=g(Z,p,p)=0$.  A
local piece of the blowup is then given by introducing a new variable
$t$ subject to the relation $f=tg$.  In this patch, $E$ is defined by
$g(Z,p,p)$.

Over $(Z,p)\in\cM\times X$, (\ref{compos}) is given by
\[
f\mapsto f(Z,p,p),\ g\mapsto g(Z,p,p),
\] 
which pulls back in the blowup to
\[
f\mapsto tg(Z,p,p),\ g\mapsto g(Z,p,p).
\] 
Dividing by the local equation of $E$ (\ref{univ}) is given by
\[
f\mapsto t,\ g\mapsto 1,
\]
which is clearly surjective as claimed.

Let $\tilde{\cI}$ be the kernel of
$\phi$.  Clearly $\tilde{\cI}$ is flat since both terms in (\ref{univ}) are,
and the fibers of $\tilde{\cI}$ over $\mathrm{Bl}_{\cC}(\cM\times X)$
are ideal sheaves in $I_{2-g}(X,\b)$. We claim that
$\tilde{\cI}\subset\cO_{\mathrm{Bl}_{\cC}(\cM\times X)\times X}$ is the
universal ideal sheaf.  

To check this, take an arbitrary family $\cJ\subset \cO_{S\times X}$
of ideal sheaves parameterized by a scheme $S$.  The saturation $\cJ'$ of 
$\cJ$ is a family of ideals in $I_{1-g}(X,\b)$ and
fits into a short exact sequence
\begin{equation}
\label{sat}
0\to \cJ  \to \cJ'\stackrel{\alpha}{\longrightarrow} \cO_{\Gamma}\to 0,
\end{equation}
where $\Gamma$ is the graph of some morphism $k:S\to X$.  The family of ideals
$\cJ'$ defines a morphism $j:S\to\cM$, leading to a morphism
\[
\psi:S\stackrel{j\times k}{\longrightarrow} \cM\times X
\]
It remains to show that $\psi$ factors through the blowup, inducing
$\cJ$ by pulling back $\tilde{\cI}$.  We must first show that the
pullback of the ideal of $\cC$ to $S$ is invertible.  We again use
local coordinates, letting $f=f(s,x),\ g=g(s,x)$ generate $\cJ'$
locally on $S\times X$, with $s\in S,\ x\in X$.  The pullback of the
ideal $I_{\cC}$ of $\cC$ is generated by $f(s,k(s))$ and $g(s,k(s))$.
Fix $s\in S$.  Without loss of generality, we can assume that the
kernel of $\alpha$ over $s$ is generated by $f$, i.e.\
$(\alpha(f))(s)=0$.  Then $f(s,k(s))=0$.  It follows that the pullback
of $I_{\cC}$ to $S$ is locally defined by the single equation
$g(s,k(s))=0$.  Hence $j\times k$ factors through the blowup, giving
the required map $S\to \mathrm{Bl}_{\cC}(\cM\times X)$.  Then the
construction of $\tilde{\cI}$ combined with (\ref{sat}) shows that
$\tilde{\cI}$ pulls back to $\cJ$ as desired, QED.

\bigskip\noindent
{\em Proof of Proposition~\ref{kkvprop}.\/}  We have
\[
I_{1-g}(X,\b)=\cM
\]
which implies
\[
D_\b^{1-g}=(-1)^{\dim \cM}e(\cM).
\]
Next 
\[
I_{2-g}(X,\b)=\mathrm{Bl}_{\cC}(\cM\times X)
\]
which implies
\[
D_\b^{2-g}=(-1)^{\dim \cM+1}\left(e(\cM)e(X)+e(\cC)\right)
\]
by a simple topological argument analogous to the argument at the end of
Section~\ref{cy}, noting that the exceptional divisor of the blowup is
a $\bP^1$ bundle over $\cC$.
Thus the $q$-expansion of the DT partition function begins
\[
Z^{DT}_\b=(-1)^{\dim \cM}e(\cM)q^{1-g}+(-1)^{\dim
\cM+1}\left(e(\cM)e(X)+e(\cC)\right)q^{2-g}+O(q^{3-g}).
\]
This leads immediately to
\[
{Z^{DT}_\b}'=Z^{DT}_\b/Z^{DT}_0=(-1)^{\dim\cM}e(\cM)q^{1-g}+
(-1)^{\dim\cM+1}e(\cC)q^{2-g}+O(q^{3-g}),
\]
QED.

\bigskip
We close by explaining how an ad hoc argument in \cite{kkv} becomes natural
and even obvious in Donaldson-Thomas theory.  A simple example will suffice.

Consider the computation of $n^0_4=-192$ in local $\bP^2$ from \cite{kkv}.  
Naive application of the formula of Conjecture~\ref{kkv} gives the incorrect
answer $n^0_4=-222$.  The explanation for the discrepancy of $-222-(-192)=-30$ 
was explained in \cite{kkv} as arising from quartic curves which are unions
of lines and cubics, i.e.\ $n^0_1n^1_3=(3)(-10)=-30$.

In Donaldson-Thomas theory, we would first compute the
Donaldson-Thomas invariants and then use (\ref{gvdt}) to solve for the 
GV invariants.  This illustrates the recursive procedure alluded to above.
{}From (\ref{gvdtexp}) and the GV invariants in \cite{kkv},
we get factors for curves of degree $d\le3$:

\begin{equation}
\label{p2dt}
\begin{array}{|c|l|}\hline
d&{\rm \ \ DT\ factor}\\ \hline
1&1+3qt+\ldots\\ \hline
2&1-6qt^2+\ldots\\ \hline
3&1-10t^3+\ldots\\ \hline
\end{array}
\end{equation}
In multiplying together the factors (\ref{gvdtexp}) for $d\le3$ we get
\begin{equation}
\label{correction}
1-30qt^4+\ldots
\end{equation}
where the omitted terms are irrelevant to the computation (although
not necessarily of higher order).  The key point is that the $d=4$
factor (\ref{gvdtexp}) contains the term $n^0_4qt^4$.  So in computing
the GV invariant, we divide $Z_{\mathrm{DT}}(X)'$ by
(\ref{correction}) and equate the coefficient of $qt^4$ with that of
the $d=4$ factor, then solve for $n^0_4$.  Thus the effect of the
division is to subtract $-30$ from the answer that would have
otherwise been obtained without the contribution of curves of lower
degree.

The conclusion is that the heuristic methods of \cite{kkv} are in
principle better adapted to the computation of DT invariants rather
than GV invariants.  In the sense above, the complications in
\cite{kkv} arise from solving for the GV invariants in terms of the DT
invariants.

\bigskip\noindent 
{\bf Acknowledgements.}  First of all, I'd like to thank the organizers
of the Snowbird conference for the nice meeting.  I'd also like to thank 
J.\ Bryan, D.\ Christie, A.\
Elezi, J.\ Guffin, A.\ Klemm, and R.\ Pandharipande for helpful conversations,
A.\ Greenspoon for corrections to an earlier version,
and the Aspen Center for Physics where part of this note was written.
My research is partially supported by NSF grants DMS 02-96154, DMS 02-44412,
and NSA grant MDA904-03-1-0050.

\end{document}